\documentclass[12pt,dvipsnames]{l4dc2021} 

\usepackage{tikz}
\usepackage{pgfplots}
\usepackage{pgf}
\usepackage{xcolor}
\usetikzlibrary{fit,calc,backgrounds}
\usepgfplotslibrary{colorbrewer}

\DeclareMathOperator{\R}{\mathbb{R}}
\DeclareMathOperator{\A}{\mathcal{A}}
\DeclareMathOperator{\B}{\mathcal{B}}
\DeclareMathOperator{\C}{\mathcal{C}}
\DeclareMathOperator{\D}{\mathcal{D}}
\DeclareMathOperator{\Han}{\mathfrak{H}}

\title[Nonlinear Data-Enabled Prediction and Control]{Nonlinear Data-Enabled Prediction and Control}
\usepackage{times}

\author{%
 \Name{Yingzhao Lian} \Email{yingzhao.lian@epfl.ch}\\
 \addr Automatic Control Lab , École polytechnique fédérale de Lausanne (EPFL) , Lausanne , Switzerland
 \AND
 \Name{Colin N. Jones} \Email{colin.jones@epfl.ch}\\
 \addr Automatic Control Lab , École polytechnique fédérale de Lausanne (EPFL) , Lausanne , Switzerland
}

\begin{document}

\maketitle

\begin{abstract}%
The Willems' fundamental lemma, which characterizes linear dynamics with measured trajectories, has found successful applications in controller design and signal processing, which has driven a broad research interest in its extension to nonlinear systems. In this work, we propose to apply the fundamental lemma to a reproducing kernel Hilbert space in order to extend its application to a class of nonlinear systems and we show its application in prediction and in predictive control.
\end{abstract}


\section{Introduction}
In the field of data-driven control, the characterization of system dynamics based on measured data serves as the driving force enabling controller design. Beyond running through a parametric modeling procedure~\citep{ljung1999system,lanzetti2019recurrent,chiuso2019system,limanond1998neural}, non-parametric methods distinguish themselves by directly representing system dynamics with data~\citep{kocijan2016modelling,hou1994parameter,guardabassi2000virtual,levine2018reinforcement}. In particular, there is a recent spark of interest in behavioral theory~\citep{willems1997introduction} where system dynamics are characterized by trajectories. This viewpoint concludes a simple and closed representation for linear time invariant systems~\citep{willems2005note}, coined Willems' fundamental lemma, and has been successfully applied to predictive control~\citep{coulson2019data}, named data-enabled predictive control. Its extension to nonlinear systems has received broad research attention, such as~\citep{rueda2020data,guo2020data,bisoffi2020data,berberich2020trajectory}. However, most works assume an explicit knowledge of the model structure with fully measured states, which in general negates the necessity of a data-driven method.

Meanwhile, as control theory is well-developed for linear systems, several works have tried to map a nonlinear control problem into a linear one. Among those trials, Koopman operator based methods establish a special viewpoint by looking into the function evolution governed by system dynamics, where the Koopman operator is a linear composite operator~\citep{koopman1932dynamical} under an autonomous system. Even though this operator is only well-defined in forward-complete systems~\citep{bittracher2015pseudogenerators}, it still serves as a successful heuristic motivating applications in system identification and controller design~\citep{lian2019learning,korda2018linear}. 

This work is inspired by both the behavioral theory and the function space viewpoints of Koopman operator theory. In particular, the linear dynamics of a linear functional is studied under the framework of the Willems' fundamental lemma, which ends up with a kernelized characterization of system responses. The contributions of this work are summarized in the following:
\begin{itemize}
    \item Propose a kernelized characterization of linear dynamics evolving in the dual space of a reproducing kernel Hilbert space and explore its application in prediction and controller design.
    \item Discuss the existence and some relevant mathematical properties of the proposed dynamics. We further elaborate the special considerations with respect to numerical implementation.
\end{itemize}

\noindent  An extended version of this work with more details is available on \href{https://arxiv.org/abs/2102.06553}{https://arxiv.org/abs/2102.06553}.

\subsection{Notation}
$\{\cdot\}_{i=1}^T$ denotes a set of size $T$ indexed by $i$, $colspan(A)$ is the column space of matrix $A$.
\section{Background}
\subsection{Reproducing Kernel Hilbert Space}
\begin{definition}\citep{saitoh2016theory}
A reproducing kernel Hilbert space (\textbf{RKHS}) over a set $X$ is a Hilbert space of functions from $X$ to $\R$ such that for each $x\in X$, the evaluation functional $E_xg:=g(x)$ is bounded. 
\end{definition}

Given an RKHS $H$, the Riesz-representation theorem~\citep{pedersen2012analysis} guarantees that each $x\in X$ corresponds to an unique $k_x\in H$ such that $<g,k_x>_H = g(x)$, where $<\cdot,\cdot>_H$ denotes the inner product defined in $H$. The kernel function defined on $H$ is given by $K(x,y)=<k_x,k_y>_H$, and it is positive-semidefinite. In particular, the dual space $H^*$ is the space of linear functionals over $H$. More background knowledge of RKHS theory is available in~\citep{berlinet2011reproducing,steinwart2008support,saitoh2016theory}.

\subsection{The Willems' Fundamental Lemma}
\begin{definition}
A Hankel matrix of depth $L$ associated with a vector signal $\{s_i\}_{i=1}^T,\;s_i\in\R^{n_s}$ is
\begingroup\makeatletter\def\f@size{9}\check@mathfonts
\begin{align*}
    \Han_L(s):=
    \begin{bmatrix}
    s_1 & s_2&\dots&s_{T-L+1}\\
    s_2 & s_3&\dots&s_{T-L+2}\\
    \vdots &\vdots&&\vdots\\
    s_L & s_{L+1}&\dots&s_T
    \end{bmatrix}
\end{align*}
\endgroup
\end{definition}

Regarding a Hankel matrix $\Han_L(s)$, the signal sequence $\{s_i\}_{i=1}^T$ is persistently exciting of order $L$ if $\Han_L(s)$ is full row rank. The Willems' fundamental lemma utilizes the Hankel matrices to characterize the response of the following deterministic linear time invariant (LTI) system, dubbed $\mathfrak{B}(A,B,C,D)$ with $x\in\mathbb{R}^{n_x}\;,\;u\in\mathbb{R}^{n_u}$ and $y\in\mathbb{R}^{n_y}$,
\begin{equation}\label{eqn:lin_dyn}
\begin{aligned}
        x_{k+1}&=Ax_k+Bu_k\\
        y_k &= Cx_k+Du_k
\end{aligned}\;, 
\end{equation}
whose order is denoted by $\mathfrak{O}(\mathfrak{B}(A,B,C,D))$ and all the $L$-step trajectories generated by this system is collected by $\mathfrak{B}_L(A,B,C,D)$. Given a sequence of input and output measurements $\{u_i\}_{i=1}^{T}$, $\{y_i\}_{i=1}^T$ , we define the following stacked Hankel matrix 
\begin{align*}
    \Han_L(u,y):=\begin{bmatrix}
    \Han_{L}(u)\\\Han_L(y)
    \end{bmatrix}\;,
\end{align*} with which one has the following \textbf{\emph{Willems' Fundamental Lemma}}.
\begin{lemma}\label{lem:funda}\citep[Theorem 1]{willems2005note}
Consider a controllable linear system and assume $u$ is persistently excited of order $N\geq \mathfrak{O}(\mathfrak{B}(A,B,C,D))+L$. Then $\text{colspan}(\Han_L(u,y))=\mathfrak{B}_L(A,B,C,D)$.
\end{lemma}

\section{Nonlinear Data-enabled Methods}\label{sect:kern}
We consider a linear dynamical system $\mathfrak{B}(\A,\B,\C\,D)$ evolving in RKHS as
\begin{align}\label{eqn:rkhs_dyn}
    \begin{split}
    f_{i+1} &= \A f_i+\B E_{u_i}\\
    E_{y_i} &= \C f_i+\D E_{u_i}\;,
    \end{split}
\end{align}
where $f_i$ is a real-valued linear functional in RKHS $H_x^*$, similarily $E_{u_i}$ and $E_{y_i}$ are evaluation functionals in RKHS $H_u^*$ and $H_y^*$. The kernels of $H_u$ and $H_y$ are $k_u(\cdot,\cdot)$ and $k_y(\cdot,\cdot)$ respectively. Meanwhile, the dynamics are modeled by bounded linear operators $\A: H_x^*\rightarrow H_x^*$, $\B: H_u^*\rightarrow H_x^*$, $\C: H_x^*\rightarrow H_y^*$ and $\D: H_u^*\rightarrow H_y^*$. It is noteworthy to point out that these dynamics are not necessarily infinite dimensional; more discussion is allocated to this point in Section~\ref{sect:pert_ext}.

Given a sequence of measurements $\{u_i\}_{i=1}^{T},\;\{y_i\}_{i=1}^T$, we have two sequences of evaluation functionals as $\{E_{u_i}\}_{i=1}^T,\; \{E_{y_i}\}_{i=1}^T$. The corresponding $n$-column Hankel matrices are:
\begingroup\makeatletter\def\f@size{9}\check@mathfonts
\begin{align}\label{eqn:fun_han}
     \Han_{L}(E_u):=&
    \begin{bmatrix}
    E_{u_1} & E_{u_2}&\dots&E_{u_{T-L+1}}\\
    E_{u_2} & E_{u_3}&\dots&E_{u_{T-L+2}}\\
    \vdots &\vdots&&\vdots\\
    E_{u_{L}} & E_{u_L+1}&\dots&E_{u_T}
    \end{bmatrix}\; ,\;
    \Han_{L}(E_y):=
    \begin{bmatrix}
    E_{y_1} & E_{y_2}&\dots&E_{y_{T-L+1}}\\
    E_{y_2} & E_{y_3}&\dots&E_{y_{T-L+2}}\\
    \vdots &\vdots&&\vdots\\
    E_{y_{L}} & E_{y_{L+1}}&\dots&E_{y_{T}}
    \end{bmatrix}\;.
\end{align}
\endgroup
For simplicity, we further define $v(\{u_i\}_{i=1}^L,\{y_i\}_{i=1}^L):=[E_{u_1},\dots E_{u_L},E_{y_1},\dots,E_{y_L}]^\top$. The Gram matrix of the stacked Hankel matrix $\Han_L(E_u,E_y):=[v_1,\dots,v_{n}]$ is then defined by 
\begin{align}\label{eqn:gram_entry}
\begin{split}
        K_{i,j} :&= k(v(u_i,y_i),v(u_j,y_j)) \\
        &= \sum\limits_{k=0}^{L-1}<E_{u_{i+k}},E_{u_{j+k}}>_{H_u^*}+\sum\limits_{k=0}^{L-1}<E_{y_{i+k}},E_{y_{j+k}}>_{H_y^*}\\ &\stackrel{(a)}{=}\sum\limits_{k=0}^{L-1}k_u(u_{i+k},u_{j+k})+\sum\limits_{k=0}^{L-1}k_y(y_{i+k},y_{j+k})\;,
\end{split}
\end{align}
where $(a)$ holds by the fact that the Hilbert space is self-dual~\citep{pedersen2012analysis}. The corresponding RKHS generated by $k(v(u_i,y_i),v(u_j,y_j))$ is constructed by the following product topology~\citep[Chapter 1.4]{berlinet2011reproducing},
 \[H^* := \underbrace{H_{u}^*\bigotimes\dots\bigotimes H_{u}^*}_{L \text{ times}} \bigotimes \underbrace{H_{y}^*\bigotimes\dots\bigotimes H_{y}^*}_{L \text{ times}}\;.\]

With the Fundamental Lemma~\ref{lem:funda}, we can state the following theorem
\begin{theorem}\label{thm:kernel_fun}
Consider a controllable linear system $\mathfrak{B}(\A,\B,\C,\D)$ and assume $E_u$ is persistently excited of order $N\geq \mathfrak{D}(\mathfrak{B}(\A,\B,\C,\D))+Ln_u$. A trajectory of length $L$ $\{\tilde{u}_i\}_{i=1}^{L}$ and $\{\tilde{y}_i\}_{i=1}^{L}\;$ is an element of $\mathfrak{B}_L(\A,\B,\C,\D)$ if and only if there exists $g\in\mathbb{R}^n$ such that
\begin{align}\label{eqn:kernel_fun}
    g^TKg+k(\tilde{v},\tilde{v})-2\sum\limits_{i=1}^n g_i k(\tilde{v},v_i) = 0\;,
\end{align}
where $\tilde{v}:=v(\{\tilde{u}_i\}_{i=1}^L,\{\tilde{y}_i\}_{i=1}^L)$.
\end{theorem}
\begin{proof}
The inputs and outputs sequence containing the evaluation functionals is an element of $\mathfrak{B}_L(\A,\B,\C,\D)$. Hence, by Fundamental Lemma~\ref{lem:funda}, $\tilde{v}\in colspan(\Han_L(E_u,E_y))$ and there exists $g=[g_1,\dots,g_n]^\top\in\mathbb{R}^n$ such that
\begingroup\makeatletter\def\f@size{9}\check@mathfonts
\begin{align*}
    \sum\limits_{i=1}^n g_iv_i = \tilde{v}_i&\stackrel{(a)}{\Longleftrightarrow} \lVert\sum\limits_{i=1}^n g_iv_i - \tilde{v}\rVert = 0 \\
    &\stackrel{(b)}{\Longleftrightarrow} \sum\limits_{i,j=1}^n g_i<v_i,v_j>g_j-2\sum\limits_{i=1}^i g_i<v_i,\tilde{v}>+<\tilde{v},\tilde{v}> = 0\\
    & \Longleftrightarrow g^TKg+k(\tilde{v},\tilde{v})-2\sum\limits_{i=1}^n g_i k(\tilde{v},v_i) = 0\;.
\end{align*}
\endgroup
$(a)$ holds by the uniqueness of the zero-element in a Hilbert space and $(b)$ follows equation~\eqref{eqn:gram_entry}.
\end{proof}
\begin{remark}\label{rmk:eval}
We intentionally do not state that $colspan(\Han_L(E_u,E_y))=\mathfrak{B}_L(\A,\B,\C,\D)$ in this theorem. One reason is that working with an arbitrary linear functional is not necessarily trivial in most RKHS. The most important reason is that the evaluation functional is not necessarily dense in the RKHS, while dynamics~\eqref{eqn:rkhs_dyn} only considers evaluation functionals. This concern can be fixed by referring to the quotient space $\mathfrak{B}_L(\A,\B,\C,\D)/\mathbb{R}$; more details can be found in the appendix of the extended version. 
\end{remark}
\begin{remark}
One may be concerned whether it is valid to apply the fundamental lemma in the dual space. Besides the original proof in~\citep{willems2005note}, Lemma~\ref{lem:funda} can be proven by showing that all trajectories span the range of a linear operator. In a similar way, the fundamental lemma holds in the dual space by considering the corresponding linear operator in the dual space.
\end{remark}

\subsection{Nonlinear Data-enabled Prediction}\label{sect:pred}
Given an input-output sequence of length $T_m$, $\tilde{u}_m:=\{\tilde{u}_1,\tilde{u}_2,\dots\,\tilde{u}_{T_m}\}$ and $\tilde{y}_m:=\{\tilde{y}_1,\tilde{y}_2,\dots\,\tilde{y}_{T_m}\}$, an open-loop prediction of length $T_p$ is to predict $\tilde{y}_p:=\{\tilde{y}_{T_m+1},\tilde{y}_{T_m+2},\dots\,\tilde{y}_{T_m+T_p}\}$ if a sequence of inputs $\tilde{u}_p:=\{\tilde{u}_{T_m+1},\tilde{u}_{T_m+2},\dots\,\tilde{u}_{T_m+T_p}\}$ is applied from $T_m+1$ to $T_m+T_p$. Theorem~\ref{thm:kernel_fun} indicates that the prediction problem is equivalent to the following optimization problem

\begin{align}\label{eqn:pred}
    \min\limits_{y_p,g}  g^TKg+k(\tilde{v},\tilde{v})-2\sum\limits_{i=1}^n g_i k(\tilde{v},v_i)\; ,
\end{align}
where $K$ is computed from the Hankel matrix $\Han_{T_m+T_p}(E_u,E_y)$ and $\tilde{v}:=v(\{\tilde{u}_i\}_{i=1}^{T_m+T_p},\{\tilde{y}_i\}_{i=1}^{T_m+T_p})$. Note that the prediction is achieved by an optimization problem instead of by solving the nonlinear equation~\eqref{eqn:kernel_fun} in order to better accommodate the numerical solvability of nonlinear equation~\eqref{eqn:kernel_fun}, the presence of measurement noise, the model mismatch and the infinite dimensionality. At the same time, it is noteworthy to point out that a solution to equation~\eqref{eqn:kernel_fun} is a global minimizer of problem~\eqref{eqn:pred}.
\begin{remark}
A kernel heuristic is mentioned in~\citep[Section V]{berberich2020trajectory}, which relies heavily on their presumed Hammerstein and/or Wiener system structure. In particular, their algorithm decouples the reconstruction of the predicted trajectory and the selection of the weight by assuming an inverse map from RKHS to the states, which is not valid for most RKHS. As discussed in Remark~\ref{rmk:eval}, our approach looks into a quotient space resulting in an unified prediction structure. We will show that Hammerstein and Wiener systems are special cases of the proposed framework in Section~\ref{sect:existence}
\end{remark}

\subsection{Nonlinear Data-enabled Predictive Control}
To convert the method in Section~\ref{sect:pred} into a predictive control scheme, the input sequence is optimized so that the corresponding output sequence is most desirable. Under a receding-horizon scheme, it leads to an optimistic bi-level problem~\citep[Chapter 2]{dempe2002foundations} as follows:
\begin{subequations}\label{eqn:mpc}
\begin{align}
\min\limits_{u_p,y_p,g}&\sum\limits_{i=0}^{T_p}\; l(u_{T_m+i},y_{T_m+1+i})\notag\\
\text{s.t.}&u_{T_m+i}\in\mathcal{U},\; y_{T_m+1+i}\in \mathcal{Y}\notag\\
& y_p \in \operatorname{arg min}_{\tilde{y}_p} g^TKg+k(\tilde{v},\tilde{v})-2\sum\limits_{i=1}^n g_i k(\tilde{v},v_i)\label{eqn:mpc_pred}\;,
\end{align}
\end{subequations}
where $l(\cdot,\cdot)$ is the stage cost and $\mathcal{U},\;\mathcal{Y}$ are constraints for control input and outputs.

\section{Discussion and Practical Issues}
In Section~\ref{sect:kern}, the theory and the applications have been built, more theoretical details and practical issues are elaborated in this section.

\subsection{On Existence of the Proposed Model}\label{sect:existence}
One obvious question is that whether the proposed model~\eqref{eqn:rkhs_dyn} makes any practical sense. Above all, the proposed model~\eqref{eqn:rkhs_dyn} includes standard linear systems~\eqref{eqn:lin_dyn} as a special instance.
\begin{lemma}\label{lem:eq_lin}
If $H_x,\;H_u,\;H_y$ are RKHS whose kernel are $k(x,y)=x^Ty$, then model~\eqref{eqn:rkhs_dyn} is equivalent to a standard linear system~\eqref{eqn:lin_dyn}. 
\end{lemma}
\begin{proof}
By Riesz representation theorem, for each $f_i\in H_x^*$ there exists a unique $x_{i}\in H_x$ such that $\forall\; \tilde{x}\in H_x$, $f_i(\tilde{x}) = <x_i^T,\tilde{x}>=x_i^T\tilde{x}$. The following proof finds the $x_{i+1}$ such that $<x_{i+1},\cdot>=f_{i+1}(\cdot)$. For any $\tilde{x}\in H_x$, equation~\eqref{eqn:rkhs_dyn} gives 
\begin{align*}
    f_{i+1}(\tilde{x}) =& \A f_i(\tilde{x})+\B E_{u_i}(\tilde{x}) \stackrel{(1)}{=} <\A x_i,\tilde{x}>+<\B u_i,\tilde{x}> \\
    =& <\A x_i+\B u_i, \tilde{x}> \Longrightarrow f_{i+1}(\cdot) = <\A x_i+\B u_i,\cdot>\;,
\end{align*}
where $(1)$ follows the definition of kernel function. In a similar way, $E_{y_i} = \C f_i+\D E_{u_i}$ can be reformulated. Hence, we conclude a standard linear system~\eqref{eqn:lin_dyn}.
\end{proof}
\begin{remark}
Lemma~\ref{lem:eq_lin} can be generalized to Hammerestein systems and Wiener systems. Without loss of generality, we consider a Hammerstein system, whose input nonlinearity enters the dynamics through an $N_\phi$ dimensional map $\phi(\cdot)$ as follows,
\begin{align*}
    x_{i+1}&=Ax_i+B\phi(u_i)\\
    y_i &= Cx_i+D\phi(u_i)\;.
\end{align*}
If $H_u$ is generated by the kernel $k(x,y)=<\phi(x),\phi(y)>_{\mathbb{R}^{N_\phi}}$ with $N_\phi$ the dimension of $\phi(\cdot)$ and $H_x\;,\;H_y$ are generated by a linear kernel, then the resulting dynamics of~\eqref{eqn:rkhs_dyn} is a Hammerstein system following a similar proof of Lemma~\ref{lem:eq_lin}.
\end{remark}

Beyond the standard linear system, the heuristic generalized linear systems considered in Koopman operator based control~\citep{korda2018linear} is also a subset of the proposed model. In particular, this generalized linear model has a form of
\begin{align*}
    \phi(x_{i+1}) &= A \phi(x_i)+Bu\\
    y_i &= C\phi(x_i)\;,
\end{align*}
where $\phi(\cdot)$ are some chosen/learnt basis functions. If there is some RKHS, such that $\phi(\cdot)\in H_x$, this model can be shown to be a special of case of the proposed model in a way similar to Lemma~\ref{lem:eq_lin}. In this case, $H_y$ will be a subspace of $H_x$~\citep[Chapter 1]{berlinet2011reproducing}. Note that if $\phi(\cdot)$ is a polynomial, then $H_x$ will be the space equipped with the exponential kernel $k(x,y)=e^{x^Ty}$. 

\subsection{On Persistent Excitation}\label{sect:pert_ext}
The assumption of persistent excitation of $E_u$ in Theorem~\ref{thm:kernel_fun} is defined according to the rank of its Hankel matrix $\Han_L(E_u)$. Checking rank of a matrix defined by functionals is not trivial, instead, the following procedure can simplify the rank calculation.
\begin{align*}
    \text{rank}(\Han_L(E_u)) = \text{rank}(\Han_L(E_u)^T\Han_L(E_u))=: K_u\;,\;
    (K_u)_{i,j} = \sum\limits_{k=0}^{L-1}k_u(u_{i+k},u_{j+k})\;.
\end{align*}
Hence, the condition of persistent excitation is determined by the rank of the corresponding Gram matrix. However, the condition of persistent excitation is only well defined for finite order dynamics, where $H_x,\;H_u,\;H_y$ are also finite dimensional, such as the RKHS corresponding to a linear kernel or a polynomial kernel. If, instead, the RKHS is infinite dimensional, persistent excitation is no longer guaranteed, and Theorem~\ref{thm:kernel_fun} only serves as a heuristic. In this case, more informative data is more preferable if we still consider the condition of persistent excitation. This claim follows the fact that $\text{rank}(K_u)\geq \text{trace}(K_u)$ when $E_{u_i}$ is within a unit ball of $H_u^*$\footnote{This condition holds for many stationary kernels, such as the RBF kernel.}. Therefore, the nuclear norm is the convex envelop of the rank function~\citep{fazel2001rank}. Given the fact the $k_u(\cdot,\cdot)$ is strictly positive definite~\citep[Lemma 4.55]{steinwart2008support},  \textit{i.e. }$k_u(u_i,u_i)>0\;\forall\;i\in\mathbb{Z}$. Therefore, the rank of $K_u$ is non-decreasing regarding the amount of data. In conclusion, a relaxed condition of persistent excitation is required for an infinite dimensional system, and therefore needs more in-depth exploration.

\subsection{On Choice of the Kernel}
Even though the proposed method is non-parametric, the choice of kernel still determines the final performance. Above all, due to the unique correspondence between the RKHS and the kernel function by Moore-Aronszajn theorem~\citep[Theorem 3]{berlinet2011reproducing}, the choice of kernel function reflects our a-priori knowledge. In particular, linear and polynomial kernels imply symmetric dynamics around $0$. Exponential kernel $k(x,y)=e^{x^Ty}$ is used when $H_x,\;H_u,\;H_y$ is spanned by polynomials,  because polynomials are dense in the corresponding RKHS~\citep[Chapter 4]{steinwart2008support}. Moreover, if the trajectories that are close to each other in the state space also show similar behavior, the RBF kernel $k(x,y)=e^\frac{\lVert x-y\rVert}{2}$ can be used.

Beyond the a-priori knowledge about the system, the choice of kernel function also affects the solvability of  Problem~\eqref{eqn:pred} and Problem~\eqref{eqn:mpc}. We observe that exponential kernel has relatively low numerical stability as two distant trajectories result in a large gradient. We also observed that the RBF kernel leads to pervasive local minima, which causes poor performance when the Problem~\eqref{eqn:mpc} and Problem~\eqref{eqn:pred} are solved by gradient based algorithms. Based on our experiments, the kernel function $k(x,y)=e^\frac{-\lVert x-y\rVert}{2}e^{x^Ty}$ in general has the best performance for examples tried so far.

Finally, even with all points discussed above, the choice of the kernel function is still non-trivial in general. If a dictionary of kernel functions $\{k(\cdot,\cdot)\}_{i=1}^{N_k}$ is available, the choice of the kernel can be optimized over the positive-weighted sum $k(\cdot,\cdot)=\sum_{i=1}^{N_k}\alpha_i k(\cdot,\cdot)$ via minimization of the prediction error in the training data set.

\subsection{On Stochastic Model with Measurement Noise}
When the data is contaminated by measurement noise, a kernel mean embedding~\citep{muandet2016kernel} can be used, which evaluates the kernel function with respect to its distribution as
\begin{align}\label{eqn:kme}
    \mathbb{E}_{\tilde{X}}k(\cdot,x)\;,
\end{align}
where $x\sim\tilde{X}$ relating to the distribution of the measurement noise. If the distribution of measurement noise is known, such as Gaussian distribution, the Equation~\eqref{eqn:kme} has closed explicit form. If noise is unknown, then~\eqref{eqn:kme} can be evaluated by the empirical distribution. 
\begin{remark}
Notice that the measurement noise in each column of the Hankel matrix is not i.i.d. Hence, the empirical distribution has slower convergence rate than $O(\frac{1}{N})$.
\end{remark}

\section{Numerical Results}
In this section, the proposed prediction and control scheme is validated by a nonlinearly damped pendulum and a bilinear DC motor. The prediction problem is solved with \textsc{Casadi}~\citep{andersson2019casadi} in \textsc{Matlab} and the bi-level predictive control optimization is solved by \textsc{fmincon} interfacing \textsc{Casadi}. It is noteworthy that the bi-level optimization cannot be solved by reformulating the lower problem with its KKT conditions, because the lower problem is non-convex. The experiments are carried out with an Intel i7-4500U(1.8 GHz) and a 1333MHz 8 GB memory.

\subsection{Damped Pendulum}
\begin{figure}[!htb]
    \centering
    \begin{subfigure}[]{}
    \centering
    \begin{tikzpicture}
    \begin{axis}[xmin=0, xmax=0.7,
    ymin=-1.6, ymax=2,
    enlargelimits=false,
    clip=true,
    grid=major,
    mark size=0.3pt,
    width=0.45\textwidth,
    height=0.32\textwidth,
    legend style={at={(0.7,0.4)},nodes={scale=0.65, transform shape}},
    ylabel = $y$]
    
    \pgfplotstableread{data/pendulum_pred.dat}{\dat}
    
    \addplot+ [ultra thick, smooth,black] table [x={t}, y={meas1}] {\dat};
  \addlegendentry{measured part}
  \addplot+ [ultra thick, smooth,BurntOrange] table [x={t}, y={pred1}] {\dat};
  \addlegendentry{prediction}
  
  \addplot+ [ultra thick, smooth,TealBlue] table [x={t}, y={real1}] {\dat};
  \addlegendentry{real outputs}
    
    \end{axis}
    \end{tikzpicture}
    \end{subfigure}
    \begin{subfigure}[]{}
    \centering
    \begin{tikzpicture}
    \begin{axis}[xmin=0, xmax=0.7,
    ymin=-1.7, ymax=1.6,
    ytick distance=0.4,
    enlargelimits=false,
    clip=true,
    grid=major,
    mark size=0.3pt,
    width=0.45\textwidth,
    height=0.32\textwidth]
    
    \pgfplotstableread{data/pendulum_pred.dat}{\dat}
    
    \addplot+ [ultra thick, smooth,black] table [x={t}, y={meas2}] {\dat};
  
  \addplot+ [ultra thick, smooth,BurntOrange] table [x={t}, y={pred2}] {\dat};
  
  \addplot+ [ultra thick, smooth,TealBlue] table [x={t}, y={real2}] {\dat};
    
    \end{axis}
    \end{tikzpicture}
    \end{subfigure}
    
        \begin{subfigure}[]{}
    \centering
    \begin{tikzpicture}
    \begin{axis}[xmin=0, xmax=0.7,
    ymin=-1.8, ymax=1.5,
    enlargelimits=false,
    clip=true,
    grid=major,
    mark size=0.3pt,
    width=0.45\textwidth,
    height=0.32\textwidth,
    ylabel = $y$,xlabel = time ($s$)]
    
    \pgfplotstableread{data/pendulum_pred.dat}{\dat}
    
    \addplot+ [ultra thick, smooth,black] table [x={t}, y={meas3}] {\dat};
  
  \addplot+ [ultra thick, smooth,BurntOrange] table [x={t}, y={pred3}] {\dat};
  
  \addplot+ [ultra thick, smooth,TealBlue] table [x={t}, y={real3}] {\dat};
    
    \end{axis}
    \end{tikzpicture}
    \end{subfigure}
    \begin{subfigure}{}
    \centering
    \begin{tikzpicture}
    \begin{axis}[xmin=0, xmax=0.7,
    ymin=-1.9, ymax=1.9,
    enlargelimits=false,
    clip=true,
    mark size=0.3pt,
    grid=major,
    width=0.45\textwidth,
    height=0.32\textwidth,xlabel= time ($s$)]
    
    \pgfplotstableread{data/pendulum_pred.dat}{\dat}
    
    \addplot+ [ultra thick, smooth,black] table [x={t}, y={meas4}] {\dat};
  
  \addplot+ [ultra thick, smooth,BurntOrange] table [x={t}, y={pred4}] {\dat};
  
  \addplot+ [ultra thick, smooth,TealBlue] table [x={t}, y={real4}] {\dat};
    
    \end{axis}
    \end{tikzpicture}
    \end{subfigure}

    \caption{Open-loop prediction of damped pendulum. An asymmetric oscillation is observed and is learnt by the proposed method in subfigure (b). Each subplot is an open-loop prediction evaluated on different data. By referring to Section~\ref{sect:pred}, the black curves are $\tilde{y}_m$ and orange curves are the open-loop prediction solved by problem~\eqref{eqn:pred}}
    \label{fig:pred_pend}
\end{figure}
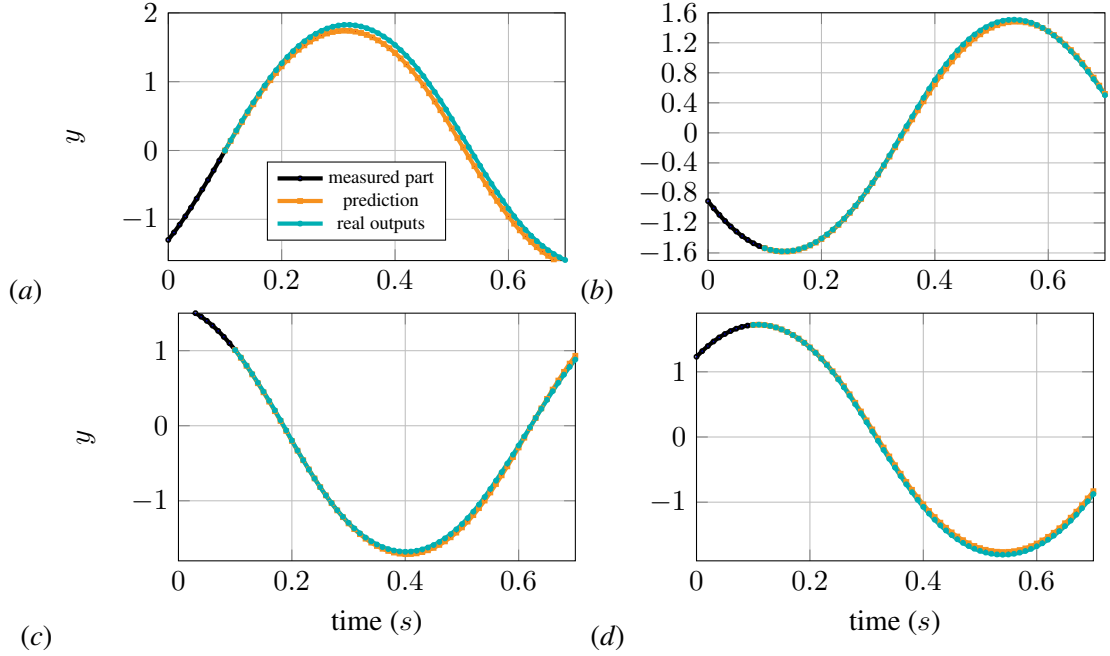

We consider a force acting on the tip of the damped pendulum. The dynamics are
\begin{align*}
    \dot{x}_1 & = x_2\;, & \dot{x}_2 &= -\frac{2g}{l}\sin(x_1)-\mu x_2^3+\frac{1}{l}|\cos(x_1)|u
\end{align*}
with $g = 9.8N/kg,\;l = 0.5m$ and $\mu = 0.1$ denoting a nonlinear friction factor. Only $x_1$ is observed as $y=x_1$. The output training data is generated by random input of ranging from -1 to 1, the sequence is measured with a sampling time of 0.04 seconds. 500 data points are used to defined the Problem~\eqref{eqn:pred} and Problem~\eqref{eqn:mpc}. Part of the output training data is shown in Figure~\ref{fig:pen_dat}, in which one can see the nonlinear modulation effect of the damped pendulum model. An open-loop prediction of $T_p=60$ step is carried out with $T_m=10$ previous step measured and the result is shown in Figure~\ref{fig:pred_pend}, where the kernel for inputs is $k_u(x,y)=0.2e^{\frac{-\lVert x-y \rVert^2}{6}}+e^{x^Ty}+0.01e^{\frac{-\lVert x-y \rVert^2}{6}}e^{x^Ty}$ and the one for outputs is $k_y(x,y) = 0.2e^{\frac{-\lVert x-y \rVert^2}{6}}+e^{x^Ty}+0.01e^{\frac{-\lVert x-y \rVert^2}{6}}e^{x^Ty}+(1+x^Ty)^2$. 

\begin{figure}
    \centering
    \begin{tikzpicture}
    \begin{axis}[xmin=0, xmax=20,
    ymin=-1.8, ymax=1.8,
    enlargelimits=false,
    clip=true,
    grid=major,
    mark size=0.3pt,
    width=0.6\textwidth,
    height=0.35\textwidth,
    ylabel= $y$,xlabel= time ($s$)]
    \pgfplotstableread{data/pendulum_data.dat}{\dat}
    \addplot+ [thick, smooth] table [x={t}, y={y}] {\dat};
    \end{axis}
    \end{tikzpicture}
    \caption{Snapshot of training data use in damped pendulum model}
    \label{fig:pen_dat}
\end{figure}
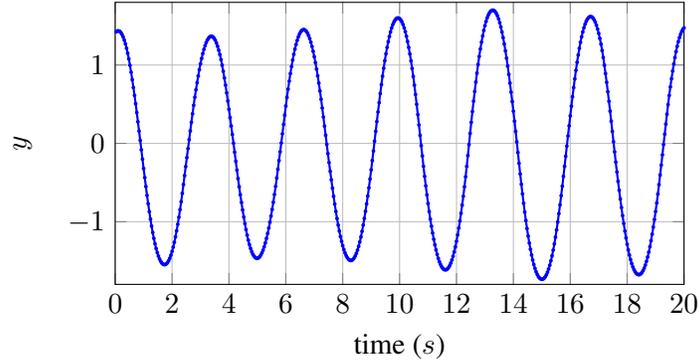

\subsection{Bilinear Motor}
\begin{figure}[!htb]
    \centering
    \begin{subfigure}[]{}
    \centering
    \begin{tikzpicture}
    \begin{axis}[xmin=0, xmax=0.7,
    ymin=-1.6, ymax= 0,
    enlargelimits=false,
    clip=true,
    grid=major,
    mark size=0.3pt,
    width=0.45\textwidth,
    height=0.3\textwidth,
    legend style={at={(0.6,0.4)},nodes={scale=0.65, transform shape}},ylabel=$y$]
    
    \pgfplotstableread{data/motor_pred.dat}{\dat}
    
    \addplot+ [ultra thick, smooth,black] table [x={t}, y={meas1}] {\dat};
    \addlegendentry{measured part}
  
  \addplot+ [ultra thick, smooth,color=BurntOrange] table [x={t}, y={pred1}] {\dat};
  \addlegendentry{prediction}
  \addplot+ [ultra thick, smooth,TealBlue] table [x={t}, y={real1}] {\dat};
    \addlegendentry{real output}
    \end{axis}
    \end{tikzpicture}
    \end{subfigure}
    \begin{subfigure}{}
    \centering
    \begin{tikzpicture}
    \begin{axis}[xmin=0, xmax=0.7,
    ymin=-1.6, ymax=-0.8,
    enlargelimits=false,
    clip=true,
    grid=major,
    mark size=0.3pt,
    width=0.45\textwidth,
    height=0.3\textwidth]
    
    \pgfplotstableread{data/motor_pred.dat}{\dat}
    
    \addplot+ [ultra thick, smooth,black] table [x={t}, y={meas2}] {\dat};
  
  \addplot+ [ultra thick, smooth,BurntOrange] table [x={t}, y={pred2}] {\dat};
  
  \addplot+ [ultra thick, smooth,TealBlue] table [x={t}, y={real2}] {\dat};
    
    \end{axis}
    \end{tikzpicture}
    \end{subfigure}
    
        \begin{subfigure}[]{}
    \centering
    \begin{tikzpicture}
    \begin{axis}[xmin=0, xmax=0.7,
    ymin=-1.3, ymax=0.6,
    enlargelimits=false,
    clip=true,
    grid=major,
    mark size=0.3pt,
    width=0.45\textwidth,
    height=0.3\textwidth,
    ylabel = $y$,xlabel = time ($s$)]
    
    \pgfplotstableread{data/motor_pred.dat}{\dat}
    
    \addplot+ [ultra thick, smooth,black] table [x={t}, y={meas3}] {\dat};
  
  \addplot+ [ultra thick, smooth,BurntOrange] table [x={t}, y={pred3}] {\dat};
  
  \addplot+ [ultra thick, smooth,TealBlue] table [x={t}, y={real3}] {\dat};
    
    \end{axis}
    \end{tikzpicture}
    \end{subfigure}
    \begin{subfigure}{}
    \centering
    \begin{tikzpicture}
    \begin{axis}[xmin=0, xmax=0.7,
    ymin=-1.4, ymax= -0.2,
    enlargelimits=false,
    clip=true,
    mark size=0.3pt,
    grid=major,
    width=0.45\textwidth,
    height=0.3\textwidth,
    xlabel = time ($s$)]
    
    \pgfplotstableread{data/motor_pred.dat}{\dat}
    
    \addplot+ [ultra thick, smooth,black] table [x={t}, y={meas4}] {\dat};
  
  \addplot+ [ultra thick, smooth, BurntOrange] table [x={t}, y={pred4}] {\dat};
  
  \addplot+ [ultra thick, smooth,TealBlue] table [x={t}, y={real4}] {\dat};
    
    \end{axis}
    \end{tikzpicture}
    \end{subfigure}

    \caption{Open-loop output prediction with random biased input sequence of bilinear motor model. Interpretation of each subplot is similar to Figure~\ref{fig:pred_pend}}
    \label{fig:motor_pred}
\end{figure}
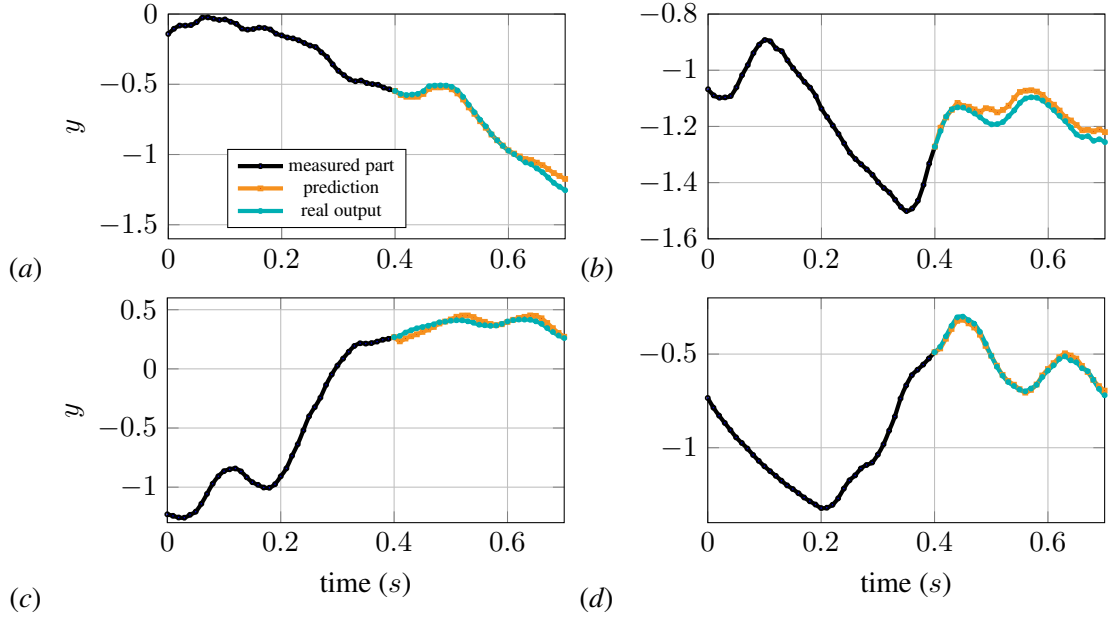

We consider a bilinear motor~\citep{daniel1998experimental} whose dynamics is
\begin{align*}
    \dot{x}_1=-\frac{R_a}{L_a}x_1 + \frac{k_m}{L_a}x_2u +\frac{u_a}{L_a}\;,\;\;\;\dot{x}_2=-\frac{B}{J}x_2 +\frac{k_m}{J}x_1u-\frac{\tau}{J}
\end{align*}
where $x_1$ is the rotor current, $x_2$ is the angular velocity and the control input $u$ is the stator current. Only the stator current is measured as $y=x_2$. The parameters are $L_a = 0.314,\;R_a = 12.345,\;k_m = 0.253,\;J = 0.00441,\; B = 0.00732,\; \tau= 1.47,$ and $u_a = 60$. Due to the synthetic effect of the bilinear term and the bias term, the responses of the system at different operating points are wildly different, hence the data is generated by a sequence of $\mathcal{N}(\mu_n,\sigma_n^2)$, where the mean of the random signal $\mu_n$ is time varying in order to excite more modes around different operating points. In particular, the mean $\mu_n$ ranges from -0.5 to 0.5, and the variance $\sigma_n=1$. 700 datapoints are used to defined the prediction and the optimal control problem. The sampling time of the generated sequence is $0.01$ seconds. An open-loop prediction of 30 is carried out with 40 previous step measured, the corresponding result is shown in Figure~\ref{fig:motor_pred}, where the kernel for inputs is $k_u(x,y)=0.1e^{\frac{-\lVert x-y \rVert^2}{4}}+e^{\frac{-\lVert x-y \rVert^2}{4}}e^{x^Ty}$ and the one for outputs is $k_y(x,y) = 0.1e^{\frac{-\lVert x-y \rVert^2}{4}}+e^{\frac{-\lVert x-y \rVert^2}{4}}e^{x^Ty}$.

\begin{figure}
    \centering
    \begin{tikzpicture}
    \begin{axis}[xmin=0, xmax=0.42,
    ymin=-1.1, ymax= -0,
    enlargelimits=false,
    clip=true,
    grid=major,
    mark size=0.5pt,
    width=0.6\textwidth,
    height=0.32\textwidth,ylabel = $y$,xlabel= time ($s$),
    legend style={font=\scriptsize}]

    \pgfplotstableread{data/motor_nmpc.dat}{\dat}
    \addplot+ [ultra thick,red] table [x={t}, y={y}] {\dat};
    \addlegendentry{Nonlinear MPC}
  
    \pgfplotstableread{data/model_ctrl.dat}{\dat}
    
    \addplot+ [ultra thick, smooth,TealBlue] table [x={t}, y={y}] {\dat};
    \addlegendentry{Output}
  
  \addplot+ [ultra thick,LimeGreen] table [x={t}, y={ref}] {\dat};
    \addlegendentry{Reference Signal}
    
    \end{axis}
    \end{tikzpicture}

    \caption{Closed-loop MPC control output}
    \label{fig:motor_ctrl}
\end{figure}
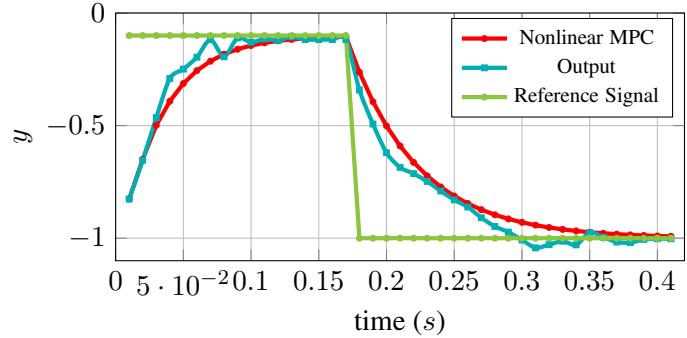

Furthermore, a predictive control scheme is tested, where the $15$ previous steps are used for a prediction horizon of $8$ steps. The stage cost is $l(u,y) = (y-y_{\text{ref}})^T(y-y_{\text{ref}})+0.01u^Tu$. The proposed predictive control scheme is compared to the nonlinear model predictive control, which has explicit knowledge of the system dynamics and full access of the state measurement. A step-like reference signal is tracked with outcome shown in Figure~\ref{fig:motor_ctrl} without considering output constraints\footnote{The consideration of output constraints makes the bi-level highly unsolvable.}, it is observed that the proposed method shows competitive performance against the model based control law, however, it fails to converge to the upper reference with a subtle bias and it has slight overshoot with respect to both set points.

\begin{remark}
The bi-level problem is currently hard to solve in real-time, it takes 15 minutes to solve the predictive control problem. Hence, a four step closed-loop control simulation takes one hour. 
\end{remark}
\section{Conclusion}
This paper presents a novel data-driven method by extending the applications of the Willems' fundamental lemma in a RKHS. The resulting kernelized characterization of system dynamics is studied and applied in prediction and predictive control with numerical validations. In conclusion, the proposed method generalizes the behavioral theory to some nonlinear dynamics, its early-stage theoretical result is presented in this work.

\section{Appendix}\label{sect:pf_thm}
As discussed in Remark~\ref{rmk:eval}, the evaluation functional is not necessarily dense in RKHS. Consider an RBF kernel, there is no $y$ such that $k(\cdot,y) = 2k(\cdot,x)$. The lack of denseness implies that not every element in $\mathfrak{B}_L(\A,\B,\C,\D)$ can be represented by a trajectory of evaluation functional, which is theoretically not desirable in behavioral theory. To fix this problem, we can consider the quotient space of $\mathfrak{B}_L(\A,\B,\C,\D)$ over $\mathbb{R}$. More specifically, it means that we consider all trajectories differs up to a non-zero scaling are considered to be equivalent. In particular the sequence $\{u_i,y_i\}_{i=0}^T$ and $\{\alpha u_i,\alpha y_i\}_{i=0}^T,\;\alpha\in\mathbb{R}$ are perceived to be equivalent in $\mathfrak{B}_L(\A,\B,\C,\D)/\mathbb{R}$ with respect to a linear dynamics~\eqref{eqn:lin_dyn}. In such space, we can proof following theorem.
\begin{theorem}
    Given an RKHS $H_x$ with kernel function $k(\cdot,\cdot)$, then the $k(\cdot,x)$ is dense in $H_x/\mathbb{R}$.
\end{theorem}
\begin{proof}
If $k(\cdot,x)$ is not dense in $H_x$, then $\exists f\in H_x$ such that there exists $\epsilon>0$ such that $\left\lVert f-\alpha k(\cdot,x)\right\rVert >\epsilon,\;\forall x\in X,\; \alpha\in\mathbb{R}$. However, the linear space span by $k(\cdot,x)$ is dense in $H_x$, hence for each $\epsilon>0$, there exist collection of $g_i$ and $x_i\in X$ such that$\left\lVert f-\sum_ig_i k(\cdot,x_i)\right\rVert<\frac{\epsilon}{2}$. Then we have
\begin{align*}
    \left\lVert f-\alpha k(\cdot,x)\right\rVert &\leq \left\lVert f-\sum_ig_i k(\cdot,x_i)\right\rVert+\left\lVert \sum_i g_i k(\cdot,x_i)-\alpha k(\cdot,x)\right\rVert\;,\\
    \Longrightarrow &\left\lVert \sum_i g_i k(\cdot,x_i)-\alpha k(\cdot,x)\right\rVert>\frac{\epsilon}{2},\;\forall \alpha\in\mathbb{R},\;\forall x\in X
\end{align*}
The last line can be lower bounded by some $\zeta>\frac{\epsilon}{2}$, such that there exist a sequence $\tilde{x}_m \in X$
\begin{align*}
    \liminf_{m\rightarrow\infty}\left\lVert \sum_i g_i k(\cdot,x_i)-\alpha k(\cdot,x)\right\rVert = \zeta\;.
\end{align*}
Then exist $M>1$ such that 
\begin{align*}
    \liminf_{m\rightarrow\infty}\frac{1}{M}\left\lVert \sum_i g_i k(\cdot,x_i)-\alpha k(\cdot,x)\right\rVert = \liminf_{m\rightarrow\infty}\left\lVert \sum_i \frac{g_i}{M} k(\cdot,x_i)-\frac{\alpha}{M} k(\cdot,x)\right\rVert<\frac{\epsilon}{2}\;.
\end{align*}

Meanwhile $\left\lVert \frac{f}{M}-\sum_i\frac{g_i}{M} k(\cdot,x_i)\right\rVert<\frac{\epsilon}{2}$ as $M>1$, in conclusion, $\liminf_{m\rightarrow\infty}\left\lVert \frac{f}{M}- \frac{1}{M} k(\cdot,x_m)\right\rVert<\epsilon$. Regarding the definition of the quotient space, we end up with a contradiction.

\end{proof}

What we discussed in this appendix only results in theoretical benefits, and has no change in the numerical implementation.

\acks{This work has received support from the Swiss National Science Foundation under the RISK project (Risk Aware Data-Driven Demand Response), grant number 200021 175627.}

\bibliography{ref}

\end{document}